\documentclass[12pt,a4paper,final]{article}
\usepackage[utf8]{inputenc}
\usepackage[english]{babel}
\usepackage{amsmath}
\usepackage{amssymb}
\usepackage{graphicx}
\usepackage[section]{placeins}
\usepackage{setspace} %use for editing only
\usepackage{authblk}
\author{Arkadiusz Kobiera}
\affil{Warsaw University of Technology}
\title{Ellipse, Hyperbola and Their Conjunction}
\date{\vspace{-5ex}}
\begin{document}
\maketitle
\begin{abstract}
This article presents a simple analysis of cones which are used to generate a given conic curve by section by a plane. It was found that if the given curve is an ellipse, then the locus of vertices of the cones is a hyperbola. The hyperbola has foci which coincidence with the ellipse vertices. Similarly, if the given curve is the hyperbola, the locus of vertex of the cones is the ellipse. In the second case, the foci of the ellipse are located in the hyperbola’s vertices. These two relationships create a kind of conjunction between the ellipse and the hyperbola which originate from the cones used for generation of these curves. The presented conjunction of the ellipse and hyperbola is a perfect example of mathematical beauty which may be shown by the use of very simple geometry. As in the past the conic curves appear to be very interesting and fruitful mathematical beings.
\end{abstract}
\section*{Introduction}

The conical curves are mathematical entities which have been known for thousands years since the first Menaechmus' research around 250 B.C. \cite{coolidge_1968}. Anybody who has attempted undergraduate course of geometry knows that ellipse, hyperbola and parabola are obtained by section of a cone by a plane. Every book dealing with the this subject has a sketch where the cone is sectioned by planes at various angles, which produces different kinds of conics. Usually authors start with the cone to produce the conic curve by section. Then, they use it to prove some facts about the conics. Many books focus on the curves themselves and their features. Even books which describe the conics theory in a quite comprehensive way \cite{coolidge_1968,otto_1976,akopyan_2007} abandon the cone after the first couple paragraphs or go to quite complex analysis of quadratics. We may find hundreds of theorems about the curves  but the relation between the cone and the conics is left to the exercise section at best \cite{otto_1976} or authors quickly go to more complex systems of conics in three-dimensional space \cite{coolidge_1968}. Probably the cone seems to be too simple to spent time on this topic, however we will show that the cone (strictly speaking family of cones) may have interesting properties as well. Apart of pure geometry, celestial mechanics is the second field where conics are important -- the orbits are conic curves. Unfortunately, the books about celestial mechanics say only a few words about the the cone if any at al. \cite{beutler_2005,wierzbinski_1973}. In this short paper we  would like to focus on the cone and its relation to conic curves which is surprisingly omitted in books, but interesting. 
\section*{Ellipse and the Cones}

\begin{figure}[h]
\centering
\includegraphics[height=10cm]{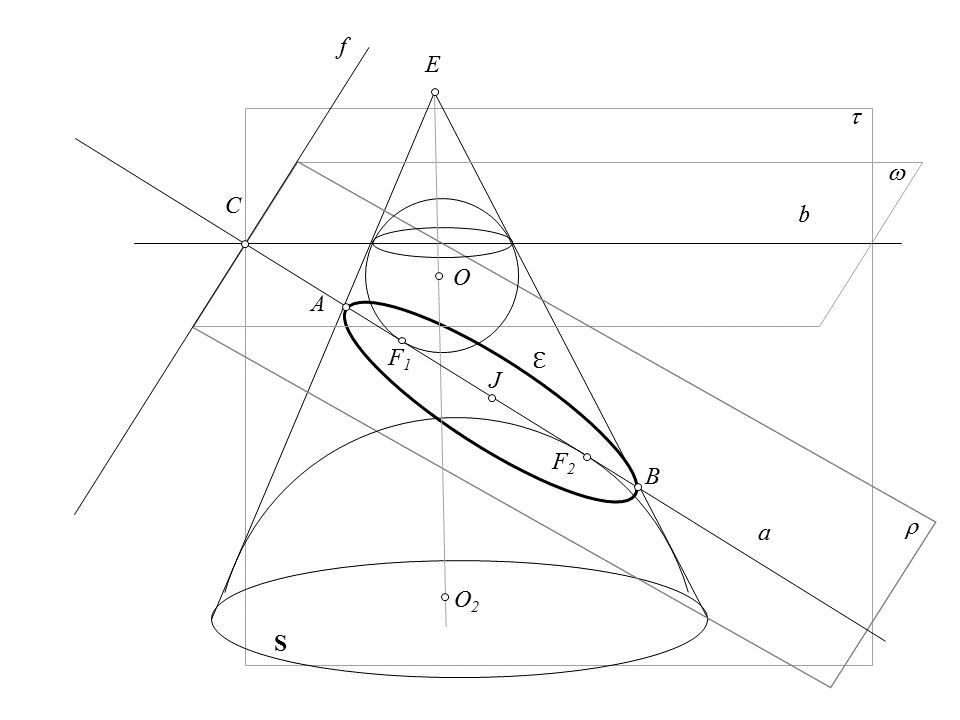}
\caption{Ellipse and the cone.}
\label{fig1}
\end{figure}

Let us consider following problem: Given is an ellipse \textbf{$\mathcal{E}$} defined by two focus points $F_{1}$ and $F_2$ and vertex $A$. This ellipse is created by section of the cone \textbf{S} by plane $\rho$. It is shown in Figure \ref{fig1}.

Our task is to find the vertex E of the cone \textbf{S}. Apart from the foci, the ellipse has also two characteristics points: the vertices $A$ and $B$. The distances from the vertices to one of the focus e.g. $F_1$ will be noted as  $r_a =|F_1A|$, and $r_b =|F_1B|$. The semi-axes of this ellipse are $a=|AJ|$ and $b=|H_1J|$ where $J$ is the center of the ellipse. The distance between foci is $c=|F_1F_2|$. The radii $r_a$ and $r_b$ define the eccentricity:
\begin{equation}
e=\frac{r_a-r_b}{r_a+r_b}=\frac{c}{a} .
\end{equation}
Obviously, we may use any set of these parameters to define the ellipse  $\mathcal{E}$, however we will prefer the radii and focus $F_1$.

\begin{figure}[h]
\centering
\includegraphics[height=10cm]{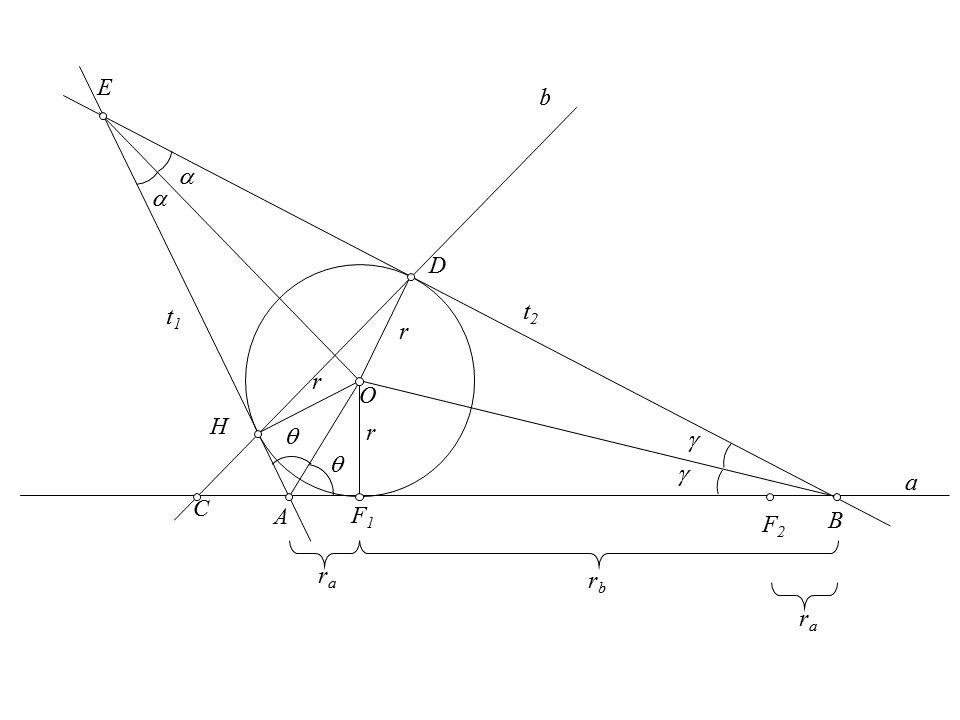}
\caption{Section of cone \textbf{S} by plane $\tau$.}
\label{fig2}
\end{figure}

The first question is about the cone: \textsl{"Is the cone \textbf{S} unique?"} The answer is in the following lemma:
\newtheorem{lem}{Lemma}
\begin{lem}
If the ellipse $\mathcal{E}$ lies on plane $\rho$ and it is defined by two vertices A, B, and focus $F_1$ (or foci $F_1$, $F_2$ and vertex A) then it may be generated by infinite number of cones \textbf{S} sectioned by the plane $\rho$.
\label{lem1}
\end{lem}
\textbf{Proof}:

The proof will be explained in a rather quite informal manner. To solve this exercise let’s reduce the three-dimensional problem to a two-dimensional problem by considering plane $\tau$ which is defined by cone's axis and foci (or vertices) of the ellipse. It is shown in Figure \ref{fig2}. We put  line \textit{a} on plane $\tau$. The line coincidences with the ellipse vertices \textit{A} and \textit{B} and the foci $F_1$, $F_2$ as well. The line $a$ is also an intersection of planes $\tau$ and $\rho$. Note that the focus points (e.g. $F_1$) are points of tangency of a sphere of center \textit{O} with the plane $\rho$. This sphere is called Dandelin’s sphere and it is simultaneously tangent to the cone. The tangency points of the sphere and the cone create a circle which defines plane $\omega$  \cite{coolidge_1968}. The intersection of planes $\omega$ and $\rho$  is line \textit{f}. We create also an additional line \textit{b} on plane $\omega$ which is perpendicular to \textit{f} and goes through the axis of the cone.  The intersection of the Dandelin's sphere by the plane $\tau$ is a circle with center \textit{O}. The circle is tangent to lines $t_1$ and $t_2$ which are two elements of the cone. These lines are obtained by cutting the cone by plane $\tau$. They meet line \textit{a} at points \textit{A} and \textit{B}. 

The problem was reduced to a problem of finding the point $E$ which is vertex of triangle \textit{ABE} that inscribes circle \textit{O}.  Since the ellipse $\mathcal{E}$ is given, the  three points \textit{A, B} and $F_1$ are fixed. Point $E$ is a point of intersection of lines $t_1$ and $t_2$. These lines are defined by points $A$ and $B$ and the circle $O$ which is tangent to the lines. If the radius \textit{r} is smaller than a certain limit $r_{max}$ then the two lines meet at point \textit{E} (this fact seems to be quite obvious so we skip this part of the proof). The maximum radius $r_{max}$ is determined by the case when the lines $t_1$ and $t_2$ are parallel. In such case the lines $t_1$ and $t_2$ become element lines of a cylinder as it is a limiting case of the cone when the point $E$ goes to infinity. In this case the ellipse $\mathcal{E}$ is obtained as a section of the cylinder. One may show that the limiting radius is equal to the minor semi-axes of the ellipse
\begin{equation}
r_{max}=a \sqrt{1-e^2}.
\end{equation}	
\noindent If the radius \textit{r} can be of any length between 0 and $r_{max}$ then the location of point \textit{E} is not unique and its position depends on radius $r$. Hence, one can construct an infinite number of cones which may be used to generate the ellipse $\mathcal{E}$.
\begin{flushright}
$\square$
\end{flushright}

\begin{figure}[h]
\centering
\includegraphics[height=10cm]{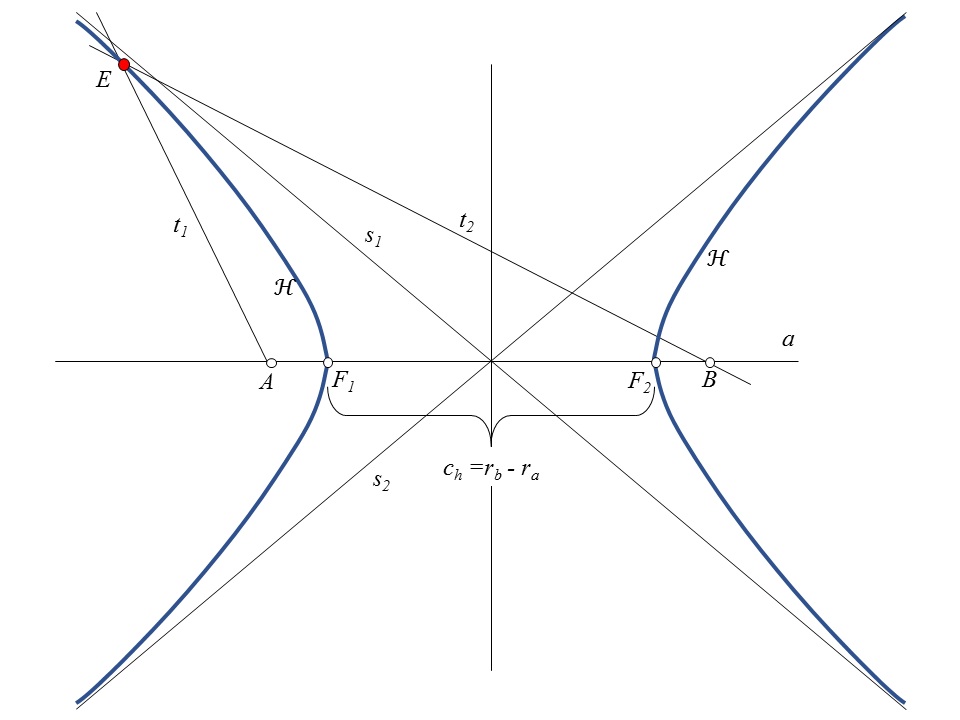}
\caption{Hyperbola $\mathcal{H}$.}
\label{fig3}
\end{figure}

If the cone S is not unique, the next question is:  \textsl{"What is the locus of the cone vertices \textit{E}?"} First, we calculate the distance from the cone vertex \textit{E} to the ellipse vertex \textit{B}

\begin{equation}
 |EB|=|BD|+|ED|=r_b+|ED|
\end{equation} 
The second equality results from the fact that \textit{BD} and $BF_1$ are tangent to circle \textit{O} and they have common endpoint \textit{B}. Obviously, the  angles $\angle F_1BO$ and $\angle DBO$ are equal and right triangles $OF_1B$ and $ODB$ are congruent. Then segments $F_1B$ and $BD$ are of the same length $r_b$. One can write similar equations for segment \textit{EA}
\begin{equation}
|EA|=|HA|+|EH|=r_a+|ED|.
\end{equation}
Here we use the fact that triangles $HOE$ and $DOE$ are congruent and triangles $HOA$ and $F_1OA$ are congruent as well. Comparison of the above equation leads to following proposition:

\newtheorem{prop}{Proposition}
\begin{prop}
\label{prop1}
If the ellipse $\mathcal{E}$ defined  by two  vertices \textit{A, B} and focus $F_1$ (or foci $F_1, F_2$ and vertex \textit{A}) is generated by section of the cones \textbf{S} by the plane $\rho$ then the locus of vertices \textit{E} of the all possible cones \textbf{S}  is a hyperbola $\mathcal{H}$. The foci of the hyperbola $\mathcal{H}$ are points \textit{A} and $B$, vertices are points $F_1$ and $F_2$ (ellipse foci). 
\end{prop}
\textbf{Prof}:

Let us calculate the difference of length of two segments \textit{EB} and \textit{EA}
\begin{equation}
|EB|-|EA|=r_b+|ED|-\left(r_a+|ED|\right)=r_b-r_a=const,
\end{equation}
\begin{equation}
|EB|-|EA|=|F_1F_2|.
\end{equation}
This difference is a constant number because $r_a$ and $r_b$ are constant as ellipse parameters, also the distance $|F_1F_2|$ is obviously constant. This directly agrees with definition of a hyperbola which foci are located in points \textit{A} and \textit{B} (see Figure \ref{fig3}). It is also clear that vertices of this hyperbola $\mathcal{H}$  are points $F_1$ and $F_2$. 
\begin{flushright}
$\square$
\end{flushright}
Indeed, the hyperbola $\mathcal{H}$ contains all the possible locations of vertices \textit{E}. The left branch contains vertices where the Dandelin's sphere is tangent to focus $F_1$. The upper part of the branch represents the case when the sphere is above the plane $\rho$, lower part is for opposite position of the sphere. The right branch is for the case where the sphere is tangent at point $F_2$. If the radius $r$ of the sphere vanishes to 0, the point $E$ goes toward foci $F_1$ or $F_2$. If the sphere’s radius $r$ goes to the maximum value $r_{max}$ the point $E$ goes to infinity on the hyperbola’s branches. Asymptotic lines $s_1$, $s_2$ are the axes of cylinders which are limiting cases of the cones with vertex in infinity. 
\section*{Hyperbola and the Cones}

\begin{figure}[h]
\centering
\includegraphics[height=10cm]{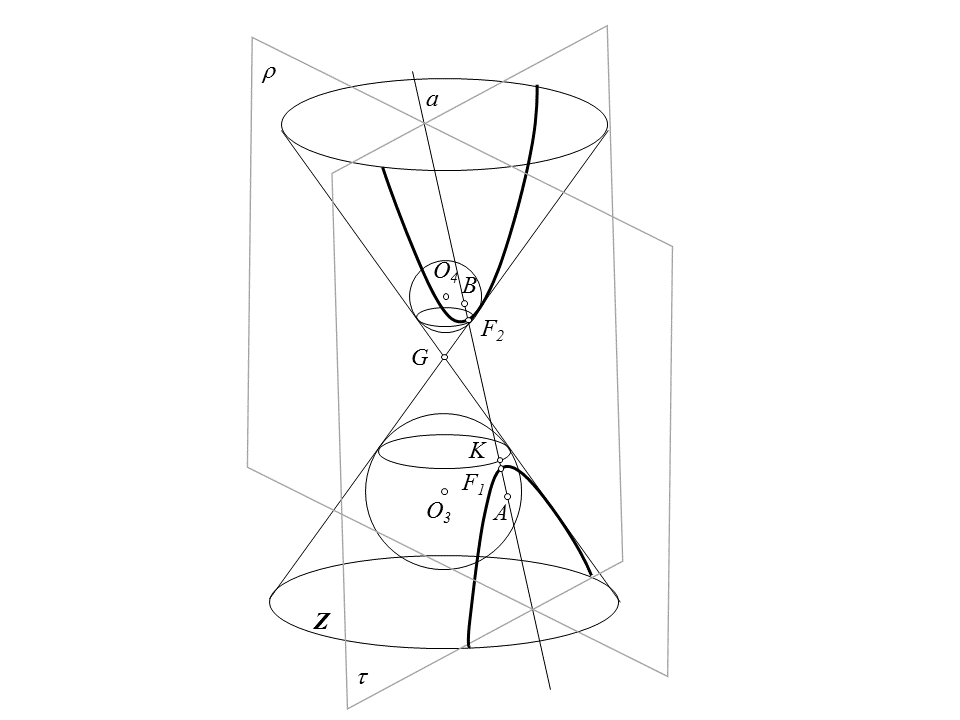}
\caption{Hyperbola $\mathcal{H}$ created by section of the cone \textbf{Z}.}
\label{fig4}
\end{figure}

Now we can ask reversed question: \textsl{What is the locus of vertices $G$ of cones \textbf{Z} which generate the given hyperbole $\mathcal{H}$.} One can consider the hyperbola $\mathcal{H}$ which was found in the previous part. This will not reduce generality of our reasoning. We will keep same plane $\tau$ where four points are defined \textit{A, B}, $F_1$ and $F_2$. They also define the hyperbola $\mathcal{H}$ on the plane $\tau$. Figure \ref{fig4} shows the situation where the hyperbola is created by sectioning the cone \textbf{Z}  by plane $\tau$.
We state the following lemma by analogy to the case of ellipse:

\begin{lem}
If the hyperbola $\mathcal{H}$ lies on plane $\tau$ and it is defined by two foci \textit{A, B,} and vertex $F_1$ (or two vertices $F_1, F_2$ and focus $A$) then it may be generated by infinite number of cones \textbf{Z} sectioned by plane $\tau$.
\label{lem2}
\end{lem}
\textbf{Proof:}

\begin{figure}[h]
\centering
\includegraphics[height=10cm]{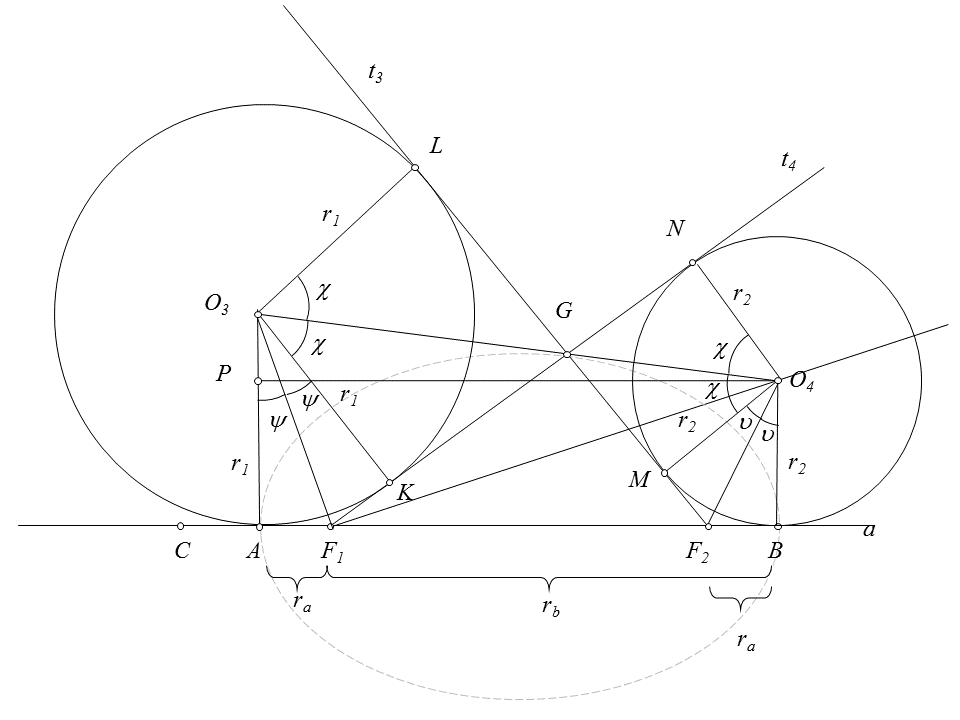}
\caption{Section of the cone \textbf{Z} by plane $\rho$.}
\label{fig5}
\end{figure}

The proof is analogous to proof of Lemma \ref{lem1}. First, we reduce the problem to planimetry by considering the plane $\rho$ (see Figure \ref{fig5}). 

By lemma's assumption we have points \textit{A, B} and $F_1$ given, then also the point $F_2$ is established because it is the vertex of the hyperbola. The vertex \textit{G} of the cone is defined by section of two lines $t_3$ and $t_4$ which are elements of cone \textbf{Z}. The lines lie on plane $\rho$ and go through points $F_1$ and $F_2$ and are tangent to two circles $O_3$ and $O_4$ respectively. The circles are sections of Dandelin's spheres (Figure \ref{fig4}) which are tangent to the plane $\tau$. They are also tangent to line $a$ at points  $A$ and $B$ which are foci of the hyperbola $\mathcal{H}$.  Let $r_1$ be the radius of the circle $O_3$. Then the point \textit{G} is not unique and its position depends on the radius $r_1$. The radius $r_1$ may vary from zero to infinity: $0 < r_1< \infty$. Hence, there exists infinite number of cones Z which generate the hyperbola if they are sectioned by the plane $\tau$.
\begin{flushright}
$\square$
\end{flushright}

The next step is finding the locus of vertices \textit{G}. By analogy to the Proposition \ref{prop1} we write following proposition:
\begin{prop}
If the hyperbola $\mathcal{H}$ defined by two vertices $F_1, F_2$ and focus \textit{A} (or foci \textit{A, B} and vertex $F_1$) is generated by section of the cones \textbf{Z} by plane $\tau$ then the locus of vertices $G$ of the cones \textbf{Z} is an ellipse $\mathcal{E}$. The foci of the ellipse $\mathcal{E}$ are points $F_1$ and $F_2$, vertices are points \textit{A} and \textit{B} (ellipse foci). 
\label{prop2}
\end{prop}
\textbf{Proof:}

We will look for the relationship between the distances from the vertex $G$ to the points $F_1$ and $F_2$. 
First, we will consider the right triangle $O_3PO_4$ (Figure \ref{fig5}). Point \textit{P} is the normal projection of point $O_4$ onto segment $AO_3$. By using Pythagoras theorem we have
\begin{equation}
|O_3 O_4 |^2=|AB|^2+(r_1-r_2)^2.
\label{eq10}
\end{equation}
where $r_2$ is radius of the circle $O_3$. The triangles $GKO_3$ and $GMO_4$ are also right triangles because the points \textit{K} and \textit{N} are points of tangency of the lines $t_3$ and $t_4$  to the circles $O_3$ and $O_4$. Hence, one can write
\begin{equation}
|O_3G|^2=|GK|^2+r_1^2,
\label{eq11}
\end{equation}
\begin{equation}
|O_4G|^2=|GM|^2+r_2^2.
\label{eq12}
\end{equation}
Recalling that
\begin{equation}
|O_3O_4 |=|O_3 G|+|O_4G|
\label{eq13}
\end{equation}
and substituting equations (\ref{eq11}), (\ref{eq12}) and (\ref{eq13}) to equation (\ref{eq10}) the following equation is obtained
\begin{multline}
(|KG|+|MG|)^2-2|KG||MG|+2|O_3G||O_4G|=|AB|^2-2r_1 r_2= \\
=|AB|^2-2|O_4 M||O_3 K|
\label{eq14}
\end{multline}
When both sides of equation (\ref{eq14}) are divided by $|O_4 M||O_3 K|$ we  get
\begin{equation}
\frac{(|GK|+|GM|)^2}{2r_1 r_2}-\left(\frac{|KG|}{|O_3 K|}\frac{|GM|}{|O_4 M|} -\frac{|O_3 G|}{|O_3 K|}\frac{|O_4 G|}{|O_4 M|}\right)=\frac{|AB|^2}{2r_1 r_2}-1.
\label{eq15}
\end{equation}
	
The triangles $O_3KG$ and $O_4MG$ are similar because they are right triangles and and angles $\angle KGO_3$ and $\angle MGO_4$ are equal. The second statement is true because the triangles  $O_4MG$ and $O_4NG$ are congruent and angles $\angle KGO_3$ and $\angle NGO_4$ are congurent as well (points \textit{G}, $O_3$ and $O_4$ lie on the axis of the cone, hence the segments $O_3G$ and $O_4G$ are co-linear). Let the measure of angles $\angle KO_3G$ and $\angle MO_4G$ be $\chi$. Simple trigonometrical relations based on Figure \ref{fig5} yield:
\begin{equation}
\frac{|KG|}{|O_3 K|}=\tan \chi=\frac{|GM|}{|O_4 M|},
\end{equation}
\begin{equation}
\frac{|O_3 G|}{|O_3 K|}=\frac{1}{\cos \chi}=\frac{|O_4 G|}{|O_4 M|}.
\end{equation}
\noindent
The second term of left hand side of equation (\ref{eq15}) can be simplified by use of the two relationship stated above
\begin{equation}
\left(\frac{|KG|}{|O_3 K|}\frac{|GM|}{|O_4 M|} -\frac{|O_3 G|}{|O_3 K|}\frac{|O_4 G|}{|O_4 M|}\right) =(\tan \chi)^2 -\frac{1}{(\cos \chi)^2}=-1.
\end{equation}

\noindent Finally, we get the simple equation:
\begin{equation}
	(|GK|+|GM|)^2=|AB|^2-4 r_1 r_2.
\end{equation}
The next step is finding the product $r_1r_2$. Let’s note that the angle $\angle A F_1 O_3$ is equal to $\angle A F_1 O_3$ and it is $\pi/2-\psi $. We may say the same about $\angle KF_1O_3$. This fact leads to conclusion that the angle $\angle NF_1B$ is equal to 
\begin{equation}
\angle NF_1B =\pi-2 \angle AF_1O_3= \pi -2(\pi /2- \psi)=2 \psi.
\end{equation}
Obviously, the line $F_1O_4$ is the bisector of this angle. Hence, the angle $\angle O_4F_1B$ is equal to $\psi$. Triangles $F_1AO_3$ and $F_1BO_4$ are similar and we may write the following proportion:
\begin{equation}
\frac{|O_3 A|}{|F_1 A|}=\frac{|F_1 B|}{|O_4 B|}.
\end{equation}
Length of segment $O_3A$ is $r_1$, lenght of segment $O_4B$ is $r_2$. We may rewrite the above equation as
\begin{equation}
\frac{r_1}{r_a}=\frac{r_b}{r_2}.
\end{equation}
Hence,  $r_1 r_2=r_a r_b$.

We successfully arrived to the conclusion that the sum of length of segments GK and GM is constant
\begin{equation}
(|GK|+|GM|)^2=|AB|^2-4r_a r_b=const. 
\end{equation}

\noindent $|AB|$ is equal to  $r_a+r_b$ then
\begin{equation}
(|GK|+|GM|)^2=(r_a+r_b )^2-4r_a r_b=(r_b-r_a )^2.
\label{eq22}
\end{equation}

The fact that $r_a = |AF_1|=|F_1K|=|F_2M|=|BF_2|$ and equation (\ref{eq22}) allows us to calculate the sum of the distances between vertex $G$ and the foci $F_1$ and $F_2$
\begin{equation}
|GF_1 |+|GF_2 |=|GK|+r_a+|GM|+r_a=(|GK|+|GM|)+2r_a=r_b+r_a.
\end{equation}
Hence, the sum of distances of the vertex $G$ from foci $F_1$ and $F_2$ is constant ($|AB|=r_a+r_b$)
\begin{equation}
|GF_1 |+|GF_2 |=|AB|.
\end{equation}
This equation is the simplest form of definition of the ellipse and we proved the proposition.
\begin{flushright}
$\square$
\end{flushright}

\section*{Conclusions}

\begin{figure}[h]
\centering
\includegraphics[height=10cm]{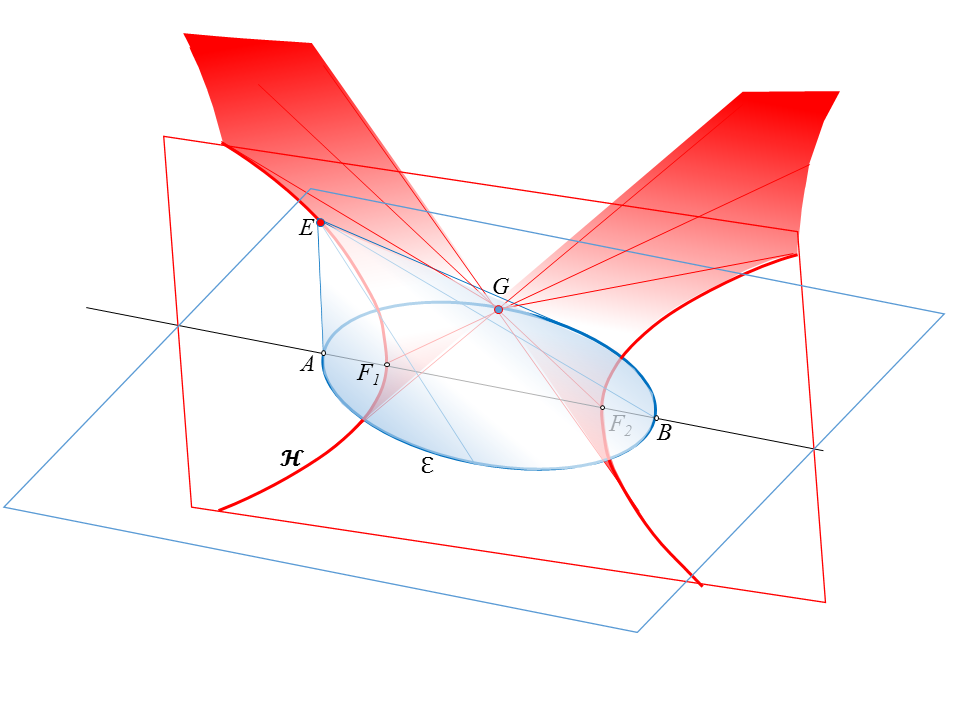}
\caption{Conjugate ellipse $\mathcal{E}$  and hyperbola $\mathcal{H}$  as curves generated by section of cones whose vertices are located on these curves.}
\label{fig6}
\end{figure}

We have shown the existence of a very interesting relationship between ellipse and hyperbola by use of very simple geometry. It was shown that ellipse and hyperbola are conjugate. This conjunction is created by locus of vertices of cones which generate the two conics. Although it seems to be a very basic property of the conics, surprisingly it is not mentioned even in some books devoted to conics geometry only. On the other hand it is wonderful that such simple mathematics may lead to such interesting results and express the beauty of geometry that is imperfectly shown in Figure \ref{fig6}.   
\bibliographystyle{plain} 
\bibliography{Kobiera-Ellipse_Hyperbola_and_Their_Conjunction}
\end{document}